\documentclass[english,a4paper,11pt]{article}
\usepackage{theorem}
\usepackage{amssymb,amsmath}
\usepackage{latexsym}
\usepackage{babel}
\usepackage{color}
\usepackage[latin1]{inputenc}
\usepackage{color,graphicx,hyperref,array}

\theoremheaderfont{\normalfont\bfseries\upshape}
\theorembodyfont{\normalfont\mdseries\slshape}
\newtheorem{theorem}{Theorem}[section]

\theorembodyfont{\normalfont\mdseries\upshape}

\newtheorem{example}[theorem]{Example}

\numberwithin{equation}{section}

\begin{document}
\pagestyle{plain} \setcounter{page}{1}
\noindent \textbf{{\LARGE Using the Generalized Collage Theorem for Estimating Unknown Parameters in Perturbed Mixed Variational Equations}}

\noindent \textbf{A.I. Garralda-Guillem\footnote{Department of Applied Mathematics, University of Granada, Granada, Spain. Email: agarral@ugr.es}, H. Kunze\footnote{Department of Mathematics and Statistics, University of Guelph, Guelph, Canada. Email: hkunze@uoguelph.ca}, D. La Torre\footnote{SKEMA Business School, Universit\'e de la C\^ote-d'Azur, Sophia Antipolis, France and Department of Economics, Management, and Quantitative Methods, University of Milan, Italy. Email: davide.latorre@skema.edu, davide.latorre@unimi.it}, M. Ruiz Gal\'an\footnote{Department of Applied Mathematics, University of Granada, Granada, Spain. Email: mruizg@ugr.es}}

\begin{abstract}
In this paper, we study a mixed variational problem subject to perturbations, where the noise term is modelled by means of  a bilinear form that has to be understood to be ``small'' in some sense. Indeed, we consider a family of such problems and provide a result that guarantees existence and uniqueness of the solution. Moreover, a stability condition for the solutions yields a Generalized Collage Theorem, which extends previous results by the same authors.
We introduce the corresponding Galerkin method and study its convergence.
We also analyze the associated inverse problem and we show how to solve it by means of the mentioned Generalized Collage Theorem and the use of adequate Schauder bases. Numerical examples show how the method works in a practical context.
\noindent
\end{abstract}


\vspace{0.5cm}

\section{Introduction: Direct vs Inverse Problem}
\label{sec1}
In applied mathematics there are always two problems associated with a mathematical model of natural phenomena, the so called direct and inverse problem. The direct problem usually refers to the determination and the analysis of the solution to a completely prescribed equation or set of equations. In many contexts, a direct problem assumes the form of differential equations subject to known
initial conditions and/or boundary conditions.
The inverse problem, instead, describes the model from the parameter estimation point-of-view. Once the model has been created and some empirical solution has been observed, it is of paramount importance to be able to determine a combination of the unknown parameters such that the induced problem admits empirical observation as an approximate solution. One can see the inverse problem as the natural opposite of a direct problem.
The study of inverse problems has attracted a lot of attention in the literature. Very often, in fact, the inverse problem is ill-posed, while the direct problem is well-posed. When a problem is well-posed, it has the properties of existence, uniqueness, and stability of the solution \cite{Hadamard}. On the other hand, an ill-posed problem lose one or more of these desirable properties. This makes the analysis of inverse problems very challenging from a numerical perspective: even when the direct problem is easily solvable, the corresponding inverse problem can be very complex and difficult to solve.

The literature is quite rich in papers proposing ad-hoc methods to address ill-posed inverse problems: These methods usually involve a minimization problem which includes a regularization term that stabilizes the numerical algorithm. One can see \cite{Keller,Kirsch,Neto,Tarantola,Tychonoff1,Tychonoff2,vogel} and the references therein to get better details about these approaches.

Quite recently other approaches have been introduced to deal with inverse problems when the corresponding direct problem can be viewed as the solution to a fixed point equation and analyzed through the well-known Banach's fixed point theorem. These approaches rely on the so-called Collage Theorem, that it is a simple consequence of the above mentioned Banach's theorem (see \cite{Banach,Ba}). In fractal imaging, these results have been used extensively to approximate a target image by the fixed point (image) of a contractive fractal transform \cite{Ba,barnsley,KuLaMeVr,Forte,LaToMe13,Lu,VrOrLaTo14}.
Over the last few years, the same philosophy has been used to deal with inverse problems for ordinary and partial differential equations. The fact that an ordinary (and even a partial) differential equation can be formulated as a fixed point equation in a specific complete metric space provides the gateway to pursuing analysis based on some of the above results.
Indeed, solution frameworks and related results have been established for case of inverse problems for different families of ordinary differential equations (see \cite{Capasso2014,kunze2,KuVr99,kunze03a,kunze03b,KuHe07}), as well as for partial differential equations (see \cite{berkunlatrui,KuLaToVr06-sub,Ku1,KuLaToLeRuGa15,Levere}).

In this paper, we explore systems of mixed variational equations, both from the direct problem and inverse problem point of view. The mixed variational formulation of a linear elliptical boundary problem is obtained from the introduction of a
new variable, usually related to any of the derivatives of the variable
original, and whose presence is justified in many cases by its applied interest. The
theoretical results, known as the \textsl{Babu\v{s}ka--Brezzi theory}, and the corresponding numerical methods, mixed finite elements, have been successfully developed in the last decades: see, for instance, \cite{bab,bof, bre,gat}.
What we discuss in this paper, instead, is a modified mixed variational problem that includes a kind of perturbation.

The paper is organized as follows. Section \ref{sec2} presents a Generalized Collage Theorem for a family of perturbed systems of mixed variational equations. Section \ref{sec3} analyzes and discusses a Galerkin numerical method for the direct problem. Section \ref{sec4} presents the formulation of the inverse problem and provides a numerical example. Section \ref{sec5} concludes the paper.

\vspace{0.5cm}

\section{Families of Mixed Variational Equations}
\label{sec2}

Unlike the classical system of mixed variational equations corresponding to the mixed variational formulation of a differential problem,  we discuss a more general version of it, which includes a certain perturbation. The perturbation term is modelled by means of a new bilinear form, that has to be interpreted to be small in some sense.
More specifically, let $E$ and $F$ be real Hilbert spaces, $a: E \times E \longrightarrow \mathbb{R}$, $b : E \times F \longrightarrow \mathbb{R}$ and $c: F \times F \longrightarrow \mathbb{R}$ be continuous bilinear forms, and $x^*:E\longrightarrow \mathbb{R}$ and $y^*:F\longrightarrow \mathbb{R}$ be linear forms. The problem under consideration is given in these terms:
\begin{equation}
\hbox{find }(x_0,y_0)\in E \times F \hbox{ such that }
\label{model}
\left\{
\begin{array}{r}
a(x_0,\cdot )+b(\cdot ,y_0)  =  x^* \\
b(x_0,\cdot)+c(y_0,\cdot)  =  y^*
\end{array}
\right..
\end{equation}
In fact, we state a more general result for a family of problems that include a stability property, \eqref{eq:bound}, which will be essential for our purposes since it will allow us to deal with a Galerkin scheme for a specific direct problem as well as with a suitable inverse problem in the next sections. Furthermore, such a stability condition, \eqref{eq:bound}, it is a \textit{Generalized Collage Theorem} that extends those in \cite{KuLaToVr06-sub} and in \cite{berkunlatrui} in the Hilbertian framework, and that in Section \ref{sec4} will be useful in order to solve an estimating parameters problem.

\bigskip

\begin{theorem}\label{th:main}
Let $J$ be a nonempty set and, for each $j \in J$, let $E_j$ and $F_j$ be real Hilbert spaces, $a_j: E_j \times E_j \longrightarrow \mathbb{R}$, $b_j: E_j \times F_j \longrightarrow \mathbb{R}$ and $c_j: F_j \times F_j \longrightarrow \mathbb{R}$ be continuous and bilinear forms, and let
\[
K_j:=\{x \in E_j : \ b_j(x,\cdot )=0\}.
\]
Suppose that
\begin{itemize}
\item[{\rm (i)}] $x \in K_j \ \wedge \ a_j(x,\cdot)_{|{K_j}}=0 \Rightarrow \ x=0$
\end{itemize}
and for some $\alpha_j,\beta_j >0$ there hold
\begin{itemize}
\item[{\rm (ii)}] $x \in K_j \ \Rightarrow \ \alpha_j \| x \| \le \|a_j(\cdot,x)_{|{K_j}} \|$,
\item[{\rm (iii)}] $y \in F \ \Rightarrow \ \beta_j \| y \| \le \| b_j(\cdot,y) \|$.
\end{itemize}
Assume in addition that
\begin{itemize}
\item[{\rm (iii)}]
\[
\rho :=\sup_{j \in J}\max \left\{ \frac{1}{\alpha_j}, \frac{1}{\beta_j} \left( 1+\frac{\|a_j\|}{\alpha_j} \right), \frac{1}{\beta_j^2}\|a_j\| \left( 1+\frac{\|a_j\|}{\alpha_j} \right) \right\} >0
\]
\end{itemize}
and that for all $j \in J$,
\begin{itemize}
\item[{\rm (iv)}] $\| c_j \| < \displaystyle \frac{1}{\rho}$.
\end{itemize}
Then, given $j \in J$ and $(x_j^*,y_j^*) \in E_j^* \times F_j^*$ there exists a unique $(x_j,y_j) \in E_j \times F_j$ such that
\begin{equation}\label{eq:system}
\left\{
\begin{array}{c}
a_j(x_j,\cdot)+b_j(\cdot,y_j)=x_j^* \\
b_j(x_j,\cdot)+c_j(y_j,\cdot)=y_j^*
\end{array}
\right..
\end{equation}
Moreover, if for each $j \in F$, $(\hat{x}_j,\hat{y}_j) \in E_j \times F_j$, then
\begin{equation}\label{eq:bound}
\inf_{j \in J}\max  \{ \|x_j-\hat{x}_j\|,\|y_j-\hat{y}_j\|\}
\le \inf_{j \in J}
\frac{\rho}{1-\rho \| c_j\|} \left( \| x_j^*-a_j(\hat{x}_j,\cdot)-b_j(\cdot ,\hat{y}_j) \| +\|y_j^* -b_j(\hat{x}_j,\cdot)-c_j(\hat{y}_j,\cdot)\| \right).
\end{equation}
\end{theorem}

\noindent {\sc Proof.} Let $j \in J$. The existence and uniqueness of solution for problem \eqref{eq:system} is a well-known fact (see, for instance \cite[Proposition 4.3.2]{bof}), but we give a sketch of the proof in order to derive also the control of the norms in \eqref{eq:bound} in a precise way.
So, let us endow the product space $E_j \times F_j$ with the norm
\[
\| (x,y) \| :=\max\{\| x \|, \| y \|\}, \qquad (x \in E_j, \ y \in F_j)
\]
and its dual space $E_j^* \times F_j^*$ with the corresponding dual norm, that is,
\[
\| (x^*,y^*) \| :=\| x^* \|+ \| y^* \|, \qquad (x^* \in E_j^*, \ y^* \in F_j^*).
\]
According to conditions (i), (ii) and (iii) and to \cite[Theorem 2.1]{gat}, the bounded and linear operator $S_j: E_j \times F_j \longrightarrow E_j^* \times F_j^*$ defined at each $(x,y) \in E_j \times F_j$ as
\[
S_j(x,y):=(a_j(x,\cdot)+b_j(\cdot,y),b_j(x,\cdot))
\]
is an isomorphism. But, in view of \cite[Theorem 2.3.5]{atk}, in order to state the existence of a unique solution for the perturbed mixed system \eqref{eq:system} it is enough to show that
\begin{equation}\label{eq:ineq}
\|S_j^{-1} \| < \frac{1}{\| c_j \|},
\end{equation}
inequality which is valid, since in view of \cite[Theorem 4.72]{gro} or \cite[Theorem 3.6]{gar-rui1} and (iv) we have that
\[
\begin{array}{rl}
\| S_j^{-1}\| & = \displaystyle \sup_{\|x^*\|+\|y^*\| \le 1}  \| S_j^{-1}(x^*,y^*)\|\\
              & \le \displaystyle \sup_{\|x^*\|+\|y^*\| \le 1}  \max \left\{ \displaystyle \frac{\| x^* \|}{\alpha_j}+\frac{1}{\beta_j}\left( 1+\frac{\| a_j \|}{\alpha_j}\right) \|y^*\|, \frac{1}{\beta_j} \left(1+ \frac{\|a_j\|}{\alpha_j} \right) \left( \|x^*\|+\frac{\|a_j\|}{\beta_j} \|y^*\|\right) \right\}  \\
              & \le \displaystyle \sup_{\|x^*\|+\|y^*\| \le 1}  \max \left\{ \displaystyle \frac{1}{\alpha_j},\frac{1}{\beta_j}\left( 1+\frac{\| a_j \|}{\alpha_j}\right), \frac{1}{\beta_j} \left(1+ \frac{\|a_j\|}{\alpha_j} \right) ,  \frac{\| a_j\|}{\beta_j^2} \left(1+ \frac{\|a_j\|}{\alpha_j} \right)  \right\}  (\|x^*\|+\|y^*\|) \\
              & \le \rho \\
              & < \displaystyle \frac{1}{\| c_j \|}.
\end{array}
\]
Furthermore, by making use of \eqref{eq:ineq} and of \cite[Theorem 4.72]{gro} or \cite[Theorem 3.6]{gar-rui1} once again, we arrive at
\begin{equation}\label{eq:aco}
\max  \{ \|x_j\|,\|y_j\|\} \le \frac{\rho}{1-\rho \| c_j\|} \left( \| x^*\| +\|y^*\| \right),
\end{equation}
where $(x_j,y_j) \in E_j \times F_j$ is the unique solution of \eqref{eq:system}.
To conclude, given $(\hat{x}_j,\hat{y}_j) \in E_j \times F_j$, since $(x_j-\hat{x}_j,y_j-\hat{y}_j)$ is the unique solution of the perturbed mixed problem
\[
\left\{
\begin{array}{c}
a_j(x_j-\hat{x}_j,\cdot)+b_j(\cdot,y_j-\hat{y}_j)=x_j^*-a_j(x,\cdot)-b_j(\cdot,\hat{y}_j) \\
b_j(x_j-\hat{x}_j,\cdot)+c_j(y_j-\hat{y}_j,\cdot)=y_j^*-b_j(\hat{x}_j,\cdot)-c_j(\hat{y}_j,\cdot)
\end{array}
\right.,
\]
then, according to inequality \eqref{eq:aco},
\[
\max  \{ \|x_j-\hat{x}_j\|,\|y_j-\hat{y}_j\|\}  \le
\frac{\rho}{1-\rho \| c_j\|} \left( \| x_j^*-a_j(\hat{x}_j,\cdot)-b_j(\cdot ,\hat{y}_j) \| +\|y_j^* -b_j(\hat{x}_j,\cdot)-c_j(\hat{y}_j,\cdot)\| \right).
\]
Finally, the arbitrariness of $j \in F$ yields \eqref{eq:bound}.
\hfill$\textcolor[rgb]{0.16,0.4,0.61}{\Box}$

\bigskip

\begin{example}\label{ex:bilap}
Given $\Omega=(0,1)^2$, $\Gamma=\partial \Omega$, $\delta \in \mathbb{R}$ and $f \in H_0^1(\Omega)$, let us consider the boundary value problem:
\begin{equation}\label{eq:varpro}
\left\{
\begin{array}{rl}
\Delta^2\psi+\delta \psi=f & \hbox{in } \Omega \\
\psi|_{\scriptscriptstyle \Gamma} =0 & \\
\Delta \psi|_{\scriptscriptstyle \Gamma} =0 &
\end{array}
\right..
\end{equation}
If one takes $w:=-\Delta \psi$, then this problem is equivalent to
\begin{equation}\label{eq:varpro2}
\left\{
\begin{array}{rl}
w+\Delta \psi=0 & \hbox{in } \Omega \\
-\Delta w+\delta \psi=f & \hbox{in } \Omega \\
\psi|_{\scriptscriptstyle \Gamma} =0 & \\
w|_{\scriptscriptstyle \Gamma} =0 &
\end{array}
\right..
\end{equation}
Then, multiplying its first equation by a test function $v \in H_0^1(\Omega)$, and integrating by part, we arrive at
\[
\int_\Omega wv- \int_\Omega \nabla w \nabla v=0.
\]
 On the other hand, when multiplying the second equation of \eqref{eq:varpro2} by a test function $\phi \in H_0^1(\Omega)$, and, proceeding as above, we write it as
\[
-\int_\Omega \nabla w \nabla \phi - \delta \int_\Omega \psi \phi =-\int_\Omega f\phi.
\]
Therefore, if we take the real Hilbert  spaces $E=F:=H_0^1(\Omega)$, the continuous bilinear forms $a: E \times E \longrightarrow  \mathbb{R}$, $b: E \times F \longrightarrow  \mathbb{R}$ and  $c: F \times F \longrightarrow  \mathbb{R}$ defined for each $w,v \in E,$ and $\phi, \psi\in F$, as
\[
a(w,v):= \int_\Omega wv,
\]
\[
b(v,\psi):=-\int_\Omega \nabla v \nabla\psi ,
\]
and
\[
c(\psi,\phi):=-\delta \int_\Omega \psi \phi,
\]
and the continuous linear forms $x^* \in E^*$ and $y^* \in F^*$ given by
\[
x^*(v):= 0 \qquad (v \in E)
\]
and
\[
y^*(\phi):=-\int_\Omega f\phi, \qquad (\phi \in F),
\]
then we have derived this variational formulation of the problem \eqref{eq:varpro}: find $(w,\psi) \in E \times F$ such that
\[
\left\{
\begin{array}{rrl}
v \in E \ \Rightarrow \ & a(w,v)+b(v,\psi) & = x^*(v) \\
w \in W \ \Rightarrow \ & b(w,\phi)+c(\psi,\phi) & =y^*(\phi)
\end{array}
\right.,
\]
which adopts the form of \eqref{eq:system} with $\mathrm{card}(J)=1$. Then, taking into account that the operator $\Delta : H_0^1 (\Omega) \longrightarrow H^{-1}(\Omega)$ is an isomorphism, it is very easy to check that, when $\delta <1$, Theorem \ref{th:main} applies and this problems admits a unique solution $(w,\psi)$ such that, for any $(\hat{w},\hat{\psi}) \in E \times F$,
\[
\max  \{ \|w-\hat{w}\|,\|\psi-\hat{\psi}\|\} \le \frac{1}{1-\delta}
\left( \| a(x,\cdot)+b(\cdot ,y) \| +\|y^* -b(x,\cdot)-c(y,\cdot)\| \right)
\]
and, in particular,
\[
\max  \{ \|w\|,\|\psi\|\} \le \frac{\|f\|}{1-\delta} .
\]
\hfill$\Box$\par

\end{example}

\bigskip

\vspace{0.5cm}

\section{The Galerkin Algorithm}\label{sec3}

Now we focus our effort on developing the Galerkin method for the perturbed mixed problem \eqref{eq:system} when $\mathrm{card}(F)=1$.

\bigskip

\begin{theorem}\label{th:discrete}

Let $E$ and $F$ be real Hilbert spaces and that $a: E \times E \longrightarrow \mathbb{R}$, $b : E \times F \longrightarrow \mathbb{R}$ and $c: F \times F \longrightarrow \mathbb{R}$ are continuous bilinear forms. Given $n \in \mathbb{N}$, let $E_n$ and $F_n$ be finite dimensional vector subspaces of $E$ and $F$, respectively, and let
\[
K_n:=\{x \in E_n: \ b(x,\cdot )_{|_{F_n}}=0\}.
\]
Let us also suppose that
\begin{itemize}
\item[{\rm (i)}] $x \in K_n \ \wedge \ a(x,\cdot)_{|{K}_n}=0 \Rightarrow \ x=0$
\end{itemize}
and there exist $\alpha_n,\beta_n >0$ such that
\begin{itemize}
\item[{\rm (ii)}] $x \in K_n \ \Rightarrow \ \alpha_n \| x \| \le \|a(\cdot,x)_{|{K}_n} \|$,
\item[{\rm (iii)}] $y \in F_n \ \Rightarrow \ \beta_n \| y \| \le \| b(\cdot,y)_{|{E}_n} \|$
\end{itemize}
and for
\begin{itemize}
\item[{\rm (iii)}]
\[
\rho_n :=\max \left\{ \frac{1}{\alpha_n}, \frac{1}{\beta_n} \left( 1+\frac{\|a\|}{\alpha_n} \right), \frac{1}{\beta_n^2}\|a\| \left( 1+\frac{\|a\|}{\alpha_n} \right) \right\} >0,
\]
\end{itemize}
there holds
\begin{itemize}
\item[(iv)] $\|c_{|{F}_n}  \| < \displaystyle \frac{1}{\rho_n}$.
\end{itemize}
Then, given $(x^*,y^*)\in E^* \times F^*$, there exists a unique $(x_n,y_n) \in E_n \times F_n$ such that
\begin{equation}\label{eqn:gal}
\left\{
\begin{array}{c}
a(x_n,\cdot)_{|{E}_n}+b(\cdot,y_n)_{|{E}_n}=x^*_{|{E}_n} \\
b(x_n,\cdot)_{|{F}_n}+c(y_n,\cdot)_{|{F}_n}=y^*_{|{F}_n}
\end{array}
\right..
\end{equation}
Furthermore, for all $(x,y) \in E \times F$ we have that
\[
\max  \{ \|x_n-x\|,\|y_n-y\|\} \le \frac{\rho_n}{1-\rho_n \| c\|} \left( \| x^*_{|{E}_n}-a(x,\cdot)_{|{E}_n}-b(\cdot ,y)_{|{E}_n} \| +\|y^*_{|{F}_n} -b(x,\cdot)_{|{E}_n} -c(y,\cdot)_{|{F}_n}\| \right).
\]

\end{theorem}

\noindent {\sc Proof.} It follows from Theorem \ref{th:main}, by means of standard arguments.
\hfill$\textcolor[rgb]{0.16,0.4,0.61}{\Box}$

\bigskip

We conclude the section by illustrating these results with the discretization of Example \ref{ex:bilap}.

\begin{example}
Let us consider the boundary value problem in Example 2.3
\begin{equation}
\left\{
\begin{array}{rl}
\Delta^2\psi+\delta \psi=f & \hbox{in } \Omega \\
\psi|_{\scriptscriptstyle \Gamma} =0 & \\
\Delta \psi|_{\scriptscriptstyle \Gamma} =0 &
\end{array}
\right.,
\end{equation}
with $\delta \in \mathbb{R}$ and $f \in H_0^1(\Omega)$. We take $\delta=1/15$, and the function
$f\in  H_0^1(\Omega)$ defined for $(x,y) \in(0, 1)^2$ in order to have the solution
$\psi_0(x,y) :=10^3(x(x-1)y(y-1))^4$.

Now let us consider the Haar system $\{ h_k\}_{k\geq 1}$ in $L^2(0,1)$, which is a Schauder basis for such real Hilbert space. Now let us build a basis for $H_0^1(0,1)$ from it: Let us define $g_1(t):=1$ and, for all $k>1$,
$$
g_k(t)=\int_0^t h_{k-1}(s) \, ds,
$$
It is easy to prove (see \cite{fu}) that the collection of function $\{g_k\}_{k\geq 1}$ is a Schauder basis for the real Hilbert space $H^1(0,1)$ and, as a consequence, $\{g0_k\}_{k\geq 1}$, where $g0_k=g_{k+2}$, is a basis for $H_0^1(0,1)$.
We now use the following bijective mapping from $\mathbb{N}$ onto $\mathbb{N}\times \mathbb{N}$ to define a bivariate  basis for $H_0^1((0,1)^2)$: let $[ \ ]$ stand for ``integer part'' and let $\sigma: \mathbb{N} \longrightarrow \mathbb{N}\times \mathbb{N}$ be the mapping given by
\begin{equation}\label{eq:sigma}
\sigma(n):=\left\{
\begin{array}{lcl}
(\sqrt{n},\sqrt{n}) & if & [\sqrt{n}]=\sqrt{n} \\
(n-[\sqrt{n}]^2,[\sqrt{n}]+1) & if & 0 <n-[\sqrt{n}]^2\leq [\sqrt{n}] \\
([\sqrt{n}]+1,n-[\sqrt{n}]^2-[\sqrt{n}]) & if & [\sqrt{n}] <n-[\sqrt{n}]^2
\end{array}
\right. .
\end{equation}
Then, the sequence $\{G0_k\}_{k\geq 1}$ defined as
$$G0_n(s,t)=g_p0(s)g_q0(t), \qquad (s,t\in (0,1))$$
where $\sigma(n)=(p,q)$, is a Schauder basis for the real Hilbert space $H_0^1((0,1)^2)$.

We can now use this basis to construct finite dimensional subspaces of the real Hilbert spaces above: For each $m\geq 1$, let us consider the finite-dimensional subspaces of $E$ and $F$
$$
E_m:=F_m:=\mathrm{span}\{G0_1,G0_2,\ldots,G0_{m}\}.
$$
Then, the corresponding discrete problem is: Find $(w_m,\psi_m)\in E_m\times F_m$, the unique solution of the discrete perturbed system
$$
\left\{ \begin{array}{lcl}
a(w_m,G0_i)+b(G0_i,\psi_m)=x^*(G0_i) & & i=1,\ldots, m,\\
b(w_m,G0_{i-m})+c(\psi_m,G0_{i-m})=y^*(G0_{i-m}) & & i=m+1,\ldots ,2m.
\end{array}\right.
$$
We show, in the following tables, the numerical results obtained  for $m=9,25,81$. The value $(w_0,\psi_0)$ denotes the exact solution of the continuous problem  with $\delta$ given above.
\begin{center}
\begin{tabular}{|c|c|c|c|}
\hline
& $m=9$ & $m=25$ & $m=81$ \\
\hline
$\|\psi_m-\psi_0\|_{L^{2}(\Omega)}$ & $1.33\times 10^{-3}$ & $9.53\times 10^{-4}$& $4.33\times 10^{-4}$\\
\hline
$\|\psi_m-\psi_0\|_{H_0^1(\Omega)}$ & $1.46\times 10^{-2}$ & $1.16\times 10^{-2}$& $7.11\times 10^{-3}$\\
\hline
$\|w_m-w_0\|_{L^{2}(\Omega)}$ & $9.41\times 10^{-2}$ & $7.09\times 10^{-2}$& $2.56\times 10^{-2}$\\
\hline
$\|w_m-w_0\|_{H_0^1(\Omega)}$ & $1.48$ & $1.22$& $7.8\times 10^{-1}$\\
\hline
\end{tabular}

\end{center}
\end{example}

\vspace{0.5cm}

\section{The Inverse Problem}
\label{sec4}

In this section we discuss the general formulation of the inverse problem for the system of mixed variational equations \eqref{model}.
Suppose that $(\hat{x}_j,\hat{y}_j) \in E_j \times F_j$ is a pair of observed/interpolated functions. Suppose, in addition, that $a_j: E_j \times E_j \longrightarrow \mathbb{R}$, $b_j: E_j \times F_j \longrightarrow \mathbb{R}$ and $c_j: F_j \times F_j \longrightarrow \mathbb{R}$ are families of bilinear forms, and $x^*_j:E_j\longrightarrow \mathbb{R}$ and $y^*_j:F_j\longrightarrow \mathbb{R}$ are families of linear forms, all them fulfilling hypotheses (i) and (ii) in Theorem \ref{th:main}.
The inverse problem can be formulated as follows: Find $\hat j\in J$, where $J$ is a compact subset of $\mathbb{R}^p$, such that $(\hat{x}_j,\hat{y}_j)$ is an approximate solution to the perturbed mixed variational system \eqref{eq:system}. Assuming that
\[
\alpha:= \inf_{j \in J} \alpha_j >0, \quad \beta:= \inf_{j \in J} \beta_j >0, \quad \delta := \sup_{j \in J} \| a_j \|, \quad \gamma:= \inf_{j \in J} \| c_j \|>0,
\]
then conditions (iii) and (iv) in Theorem \ref{th:main} are valid as soon as
$\rho \gamma <1$, and so, such a result applies.

Then, in view of the collage estimation \eqref{eq:bound}, the inverse problem can be solved by minimizing the following objective function
\begin{equation}
\xi(j) := \| x^*_j -a_j(\hat w,\cdot)-b_j(\cdot,\hat \psi) \| + \| y^*_j - b_j (\hat w, \cdot ) -c_j(\hat \psi, \cdot) \|
\label{objective}
\end{equation}
over $j\in J$. This objective function measures the distance between the left and the right hand-side of Eq. (\ref{eq:system}). The optimal value is closer to zero
the better the approximation will be as the distance between the target solution $(\hat{x}_j,\hat{y}_j)$ and the theoretical one $(x_j,y_j)$ gets very small.
The optimization problem can be discretized by means of Schauder bases in the real Hilbert spaces involved, along the lines of \cite[Section 3]{berkunlatrui} and \cite[Section 4]{KuLaToVr06-sub}, and the minimization algorithm has been implemented using the MAPLE 2018 optimization toolbox. The optimal solution provides the estimation of the unknown parameters of the model.

\begin{figure}[h]
\centering
\includegraphics[width=6cm]{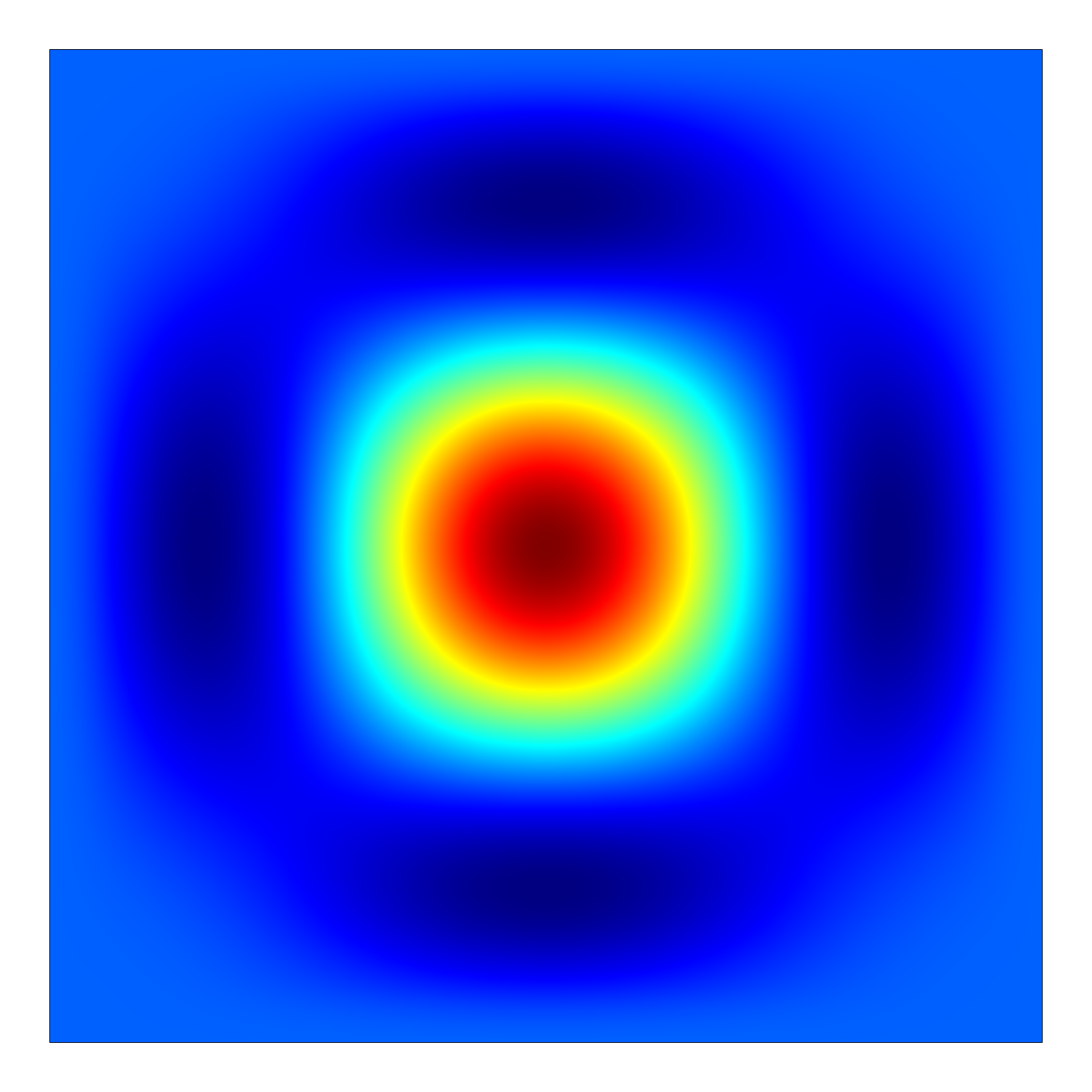}
\includegraphics[width=6cm]{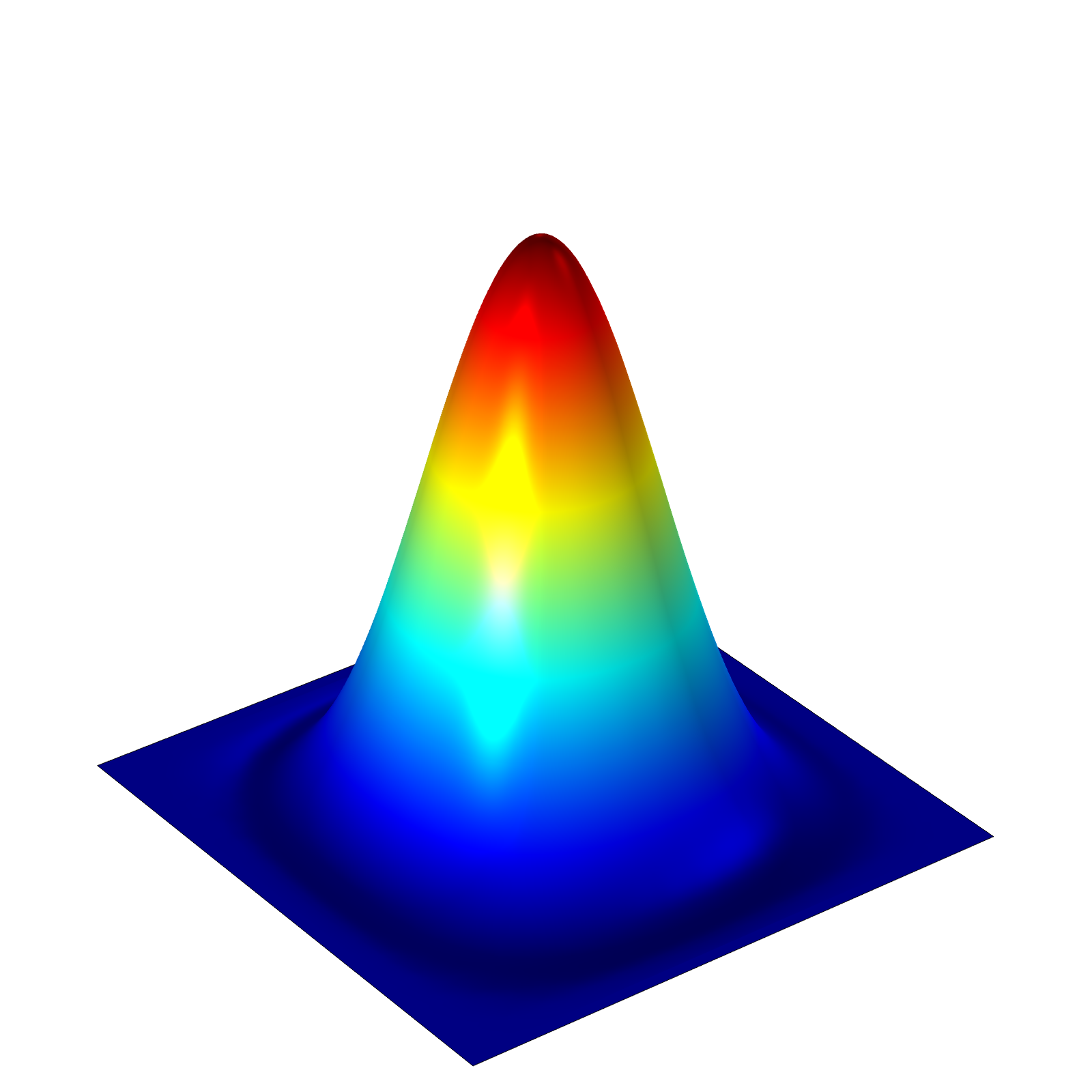}
\caption{Isotherms and surface contour plot for the target solution in the example.}
\label{fig:results}
\end{figure}

Now we illustrate a numerical implementation of the algorithm.
We start with the system in the Example \ref{ex:bilap}, setting $\delta={1\over 4}$ and choosing $f(x,y)$ such that the solution $u(x,y)$ to the problem is
$10^3 [x(1-x)y(1-y)]^4$. We solve the system in COMSOL. Isotherms and surface contour plots are shown in Figure~\ref{fig:results}.
Then we sample the numerical solution on a uniform grid of $9\times 9$ interior points of $[0,1]^2$. We interpolate each set of 81 points, with low-amplitude relative noise added, to build two target functions $\hat{u}$ and $\hat{w}$. We feed these representations into our Generalized Collage Theorem machinery; Eq. \eqref{objective} is finite dimensionalized by working with a uniform finite-element basis on $[0,1]$ with 81 interior nodes. Finally, knowing $f(x,y)$, we recover $C_1$, $C_2$, $C_3$ so that $\hat{u}$ and $\hat{w}$ are approximate solutions to the system
$$
\left\{
  \begin{array}{ll}
     C_1\Delta{u}+C_2w = 0, & \hbox{} \\
    -C_1\Delta{w}+C_3u = f(x,y), & \hbox{}
  \end{array}
\right.
$$
The true values are $C_1=1$, $C_2=1$,$C_3={1\over 4}$. The results are presented in Table~\ref{tbl:results}. The number in the final column of the table is the value of the generalized collage distance. We say that for low relative noise values, the method does reasonably well.
\begin{center}
\begin{table}
\begin{tabular}{|c|c|c|c|c|}
\hline
noise & $C_1$    & $C_2$     & $C_3$      & Collage Distance \\
\hline
0\% & 1.000107013477 & 0.9999842592864 & 0.46822577278 & 0.00054485531439904 \\
\hline
0.5\% & 1.000059971746 & 1.0000823273167 & 0.30727953435 & 0.00056207875045237 \\
\hline
1.0\% & 1.000009640383 & 1.0001771911339 & 0.15302707869 & 0.00059892386287046 \\
\hline
1.5\% & 0.999956019932 & 1.0002688495507 & 0.00547416104 & 0.00065538437959271 \\
\hline
2\% & 0.999899110999 & 1.0003573014412 & -0.13537360557 & 0.00073145365856316 \\
\hline
\end{tabular}
\caption{Results of the Numerical Simulation. True values are $(C_1,C_2,C_3)=(1,1,{1\over 4})$.}
\label{tbl:results}
\end{table}
\end{center}

One can easily notice that the estimation of the coefficient $C_3$ is not so good: this is depending on the numerical approximation of the
$\Delta u$ rather than the method itself. When solving an inverse problem, in fact, empirical data and observations for $u$ are used to estimate the unknown parameters. In this model, however, the empirical data is used to get a numerical approximation of $\Delta u$ which turns out to add more noise
to the inverse problem implementation. 

%\vspace{0.5cm}

\section{Conclusion}
\label{sec5}

In this paper we have studied the direct problem and the inverse problem for perturbed mixed variational equations. We have shown conditions that guarantee the existence and uniqueness of the solution to the direct problem and formulated the
inverse problem as an optmization problem using an extension of the Collage Theorem. We have also provided a numerical Galerkin scheme to approximate the solution to this model.
A potential application to a fourth-order PDE example is also illustrated: by substitution one can reduce this example to a perturbed mixed variational problem and then use the theory and the numerical treatment presented in this work to solve it.

\vspace{0.5cm}

\section*{Acknowledgement}
Research partially supported by project MTM2016-80676-P (AEI/FEDER, UE) and by Junta de Andaluc\'{\i}a Grant FQM359.


\vspace{0.5cm}

\begin{thebibliography}{9}

\bibitem{atk} K. Atkinson, W. Han, \textsl{Theoretical numerical analysis. A functional analysis framework}, third edition, Texts in Applied Mathematics 39, Springer, Dordrecht, 2009.

\bibitem{bab} I. Babu\v{s}ka, Error-bounds for finite element method, Numerische Mathematik 16 (1971),  322--333.


\bibitem{Banach}
Banach S, Sur les op\'{e}rations dans les ensembles abstraits et leurs applications aux \'{e}quations int\'{e}grales, Fundamenta Mathematicae 3 (1922), 133--181.

\bibitem{Ba}
Barnsley M, Fractals Everywhere, Academic Press, New York, 1989.

\bibitem{barnsley} M.F. Barnsley, V. Ervin, D. Hardin, J. Lancaster, Solution of an inverse problem for
fractals and other sets, Proceedings of the National Academy of Sciences of the United States of America 83 (1985), 1975--1977.

\bibitem{berkunlatrui} M.I. Berenguer, H. Kunze, D. La Torre, M. Ruiz Gal\'an, Galerkin method for constrained variational equations
and a collage-based approach to related inverse problems, J. Comput. Appl. Math. 292 (2016), 67--75.

\bibitem{bof} D. Boffi et al., Mixed finite elements, compatibility conditions and applications, Lecture Notes in Mathematics 1939,
Springer-Verlag, Berlin, 2008.

\bibitem{bre} F. Brezzi, On the existence, uniqueness and
approximation of saddle-point problems arising from Lagrangian
multipliers, Revue Fran\c aise d'Automatique, Informatique et Recherche Op\'erationnelle S\'erie Rouge 8 (1974),
129--151.

\bibitem{Capasso2014}
V. Capasso, H. Kunze, D. La Torre, E.R. and Vrscay, Solving inverse problems for differential equations by a ``generalized collage'' method and application to a mean field stochastic model, Nonlinear Analysis: Real World Applications 15 (2014), 276--289.

\bibitem{fu} S. Fu\v{c}ik, Fredholm alternative for nonlinear operators in Banach spaces and its applications to differential and integral equations, \v{C}asopis Pro P\v{e}stov\'an\'i Matematik  96 (1971), 371--390.

\bibitem{gar-rui1}
A.I. Garralda-Guillem, M. Ruiz Gal\'an, Mixed variational formulations in locally convex spaces, Journal of Mathematical Analysis and Applications
414 (2014), 825--849.

\bibitem{gat}
G.N. Gatica, A simple introduction to the mixed finite element method. Theory and applications,
SpringerBriefs in Mathematics, Springer, Cham, 2014.

\bibitem{gro} C. Grossmann, H.G. Roos, \textsl{Numerical treatment of partial differential equations}, Universitext, Springer, Berlin, 2007.

\bibitem{kunze2}
H. Kunze, J. Hicken, E.R. Vrscay, Inverse problems for ODEs using contraction maps: Suboptimality of the ``collage method,''
Inverse Problems 20 (2004),  977--991.

\bibitem{KuVr99}
H. Kunze, E.R. Vrscay, Solving inverse problems for ordinary differential equations using the Picard contraction mapping, Inverse Problems 15 (1999), 745--770.

\bibitem{kunze03a}
H. Kunze, S. Gomes, Solving An Inverse Problem for Urison-type Integral Equations Using Banach's Fixed Point Theorem, Inverse Problems 19 (2003), 411--418.

\bibitem{kunze03b}
H. Kunze, J. Hicken, E.R. Vrscay, Inverse Problems for ODEs Using Contraction Maps: Suboptimality of the ``Collage Method,'' Inverse Problems 20 (2004),
977--991.

\bibitem{KuHe07}
H. Kunze, K. Heidler, The Collage Coding Method and its Application to an Inverse Problem for the Lorenz System, Applied Mathematics and Computation 186 (2007), 124--129.

\bibitem{KuLaToVr06-sub}
H. Kunze, D. La Torre, E. R. Vrscay, A generalized collage method based upon the Lax--Milgram functional for solving boundary value inverse problems,  Nonlinear Analysis 71 (2009), e1337--e1343.

\bibitem{Ku1}
H. Kunze, D. La Torre, E. R. Vrscay, Solving inverse problems for variational equations using the ``generalized collage methods,'' with applications to boundary value problems, Nonlinear Analysis Real World Applications, 11 (2010), 3734--3743.

\bibitem{KuLaMeVr}
H. Kunze, D. La Torre, F.Mendivil, E. R. Vrscay, Fractal-based methods in analysis, Springer, 2012.

\bibitem{KuLaToLeRuGa15}
H. Kunze, D. La Torre, K. Levere, M. Ruiz Gal\'an, Inverse problems via the ``generalized collage theorem'' for vector-valued Lax-Milgram-based variational problems, Mathematical Problems in Engineering (2015), Art. ID 764643, 8 pp.

\bibitem{Forte}
B. Forte, E.R. Vrscay, Inverse problem methods for generalized fractal transforms. In Fractal Image Encoding and Analysis, NATO ASI Series F, Vol 159, ed. Y.Fisher, Springer Verlag, New York, 1998.

\bibitem{Hadamard}
J. Hadamard, Lectures on the Cauchy problem in linear partial differential equations, Yale University Press, 1923.

\bibitem{Keller}
J.B. Keller, Inverse Problems,  The American Mathematical Monthly 83 (1976), 107--118.

\bibitem{Kirsch}
A. Kirsch, An introduction to the mathematical theory of inverse problems, Springer, 2011.

\bibitem{LaToMe13}
D. La Torre, F. Mendivil, A Chaos Game Algorithm for Generalized Iterated Function Systems Applied Mathematics and Computation 224 (2013), 238--249.

\bibitem{Levere}
K. Levere, H. Kunze, D. La Torre, A collage-based approach to solving inverse problems for second-order nonlinear parabolic PDEs, Journal of Mathematical Analysis and Applications, 406 (2013), 120--133.

\bibitem{Lu}
N. Lu, Fractal imaging,  Morgan Kaufmann Publishers Inc., 1997.

\bibitem{Neto}
F.D. Moura Neto, A.J. da Silva Neto, An Introduction to Inverse Problems with Applications, Springer, New York, 2013.

\bibitem{Tarantola}
A. Tarantola, Inverse Problem Theory and Methods for Model Parameter Estimation, SIAM, Philadelphia, 2005.

\bibitem{Tychonoff1}
A.N. Tychonoff, Solution of incorrectly formulated problems and the regularization method, Doklady Akademii Nauk SSSR 151 (1963), 501--504.

\bibitem{Tychonoff2}
A.N. Tychonoff, N.Y. Arsenin, Solution of Ill-posed Problems, Washington: Winston \& Sons, 1977.

\bibitem{vogel}
C.R. Vogel, Computational Methods for Inverse Problems, SIAM, New York, 2002.

\bibitem{VrOrLaTo14}
E.R. Vrscay, D. Otero, D. La Torre, A Simple Class of Fractal Transforms for Hyperspectral Images, Applied Mathematics and Computation 231 (2014), 435--444.

\end{thebibliography}
\end{document}